\documentclass[aos,MSNbibl,nameyear,dvips]{arximspdf}

\doi{10.1214/12-AOS1051} 
\volume{40}
\issue{6}
\pubyear{2012}
\firstpage{2798}
\lastpage{2822}

\makeatletter
\newcommand{\rrvert}{\vert}
\newcommand{\llvert}{\vert}

\newtheorem{theorem}{Theorem}[section]
\newtheorem{lemma}{Lemma}[section]
\newtheorem{corollary}{Corollary}[section]

\newproclaim{example}{Example}[section]
\newproclaim{remark}{Remark}[section]
\newproclaim{algorithm}{Algorithm}[section]

\makeatother

\begin{document}
\begin{frontmatter}

\title{On the uniform asymptotic validity of subsampling and the bootstrap}
\runtitle{Uniform asymptotic validity}

\begin{aug}
\author[A]{\fnms{Joseph P.} \snm{Romano}\corref{}\ead[label=e1]{romano@stanford.edu}}
\and
\author[B]{\fnms{Azeem M.} \snm{Shaikh}\ead[label=e2]{amshaikh@uchicago.edu}}
\runauthor{J. P. Romano and A. M. Shaikh}
\affiliation{Stanford University and University of Chicago}
\address[A]{Departments of Economics\\
\quad and Statistics\\
Stanford University\\
Sequoia Hall\\
Stanford, California 94305-4065\\
USA\\
\printead{e1}} 
\address[B]{Department of Economics\\
University of Chicago \\
1126 E. 59th Street\\
Chicago, Illinois 60637\\
USA\\
\printead{e2}}
\end{aug}

\received{\smonth{4} \syear{2012}}
\revised{\smonth{9} \syear{2012}}

%
\begin{abstract}
This paper provides conditions under which subsampling and the
bootstrap can be used to construct estimators of the quantiles of the
distribution of a root that behave well uniformly over a large class of
distributions $\mathbf P$. These results are then applied (i) to
construct confidence regions that behave well uniformly over $\mathbf
P$ in the sense that the coverage probability tends to at least the
nominal level uniformly over $\mathbf P$ and (ii) to construct tests
that behave well uniformly over $\mathbf P$ in the sense that the size
tends to no greater than the nominal level uniformly over $\mathbf P$.
Without these stronger notions of convergence, the asymptotic
approximations to the coverage probability or size may be poor, even in
very large samples. Specific applications include the multivariate
mean, testing moment inequalities, multiple testing, the empirical
process and $U$-statistics.
\end{abstract}

%

\end{frontmatter}

\section{Introduction} \label{sectionintro}

Let $X^{(n)} = (X_1,\ldots, X_n)$ be an i.i.d. sequence of random
variables with distribution $P \in\mathbf P$, and denote by $J_n(x,P)$
the distribution of a real-valued root $R_n = R_n(X^{(n)},P)$ under
$P$. In statistics and econometrics, it is often of interest to
estimate certain quantiles of $J_n(x,P)$. Two commonly used methods for
this purpose are subsampling and the bootstrap. This paper provides
conditions under which these estimators behave well uniformly over
$\mathbf P$. More precisely, we provide conditions under which
subsampling and the bootstrap may be used to construct estimators $\hat
c_n(\alpha_1)$ of the $\alpha_1$ quantiles of $J_n(x,P)$ and $\hat
c_n(1-\alpha_2)$ of the $1-\alpha_2$ quantiles of $J_n(x,P)$, satisfying
%
\begin{equation}
\label{eqtwosided} \liminf_{n \rightarrow\infty} \inf_{P \in\mathbf{P}}
P\bigl\{\hat
c_n(\alpha_1) \leq R_n \leq\hat
c_n(1-\alpha_2)\bigr\} \geq1 - \alpha_1 -
\alpha_2.
\end{equation}
Here, $\hat c_n(0)$ is understood to be $-\infty$, and $\hat c_n(1)$
is understood to be $+\infty$.
For the construction of two-sided confidence intervals of nominal level
$1- 2 \alpha$ for a real-valued parameter,
we typically would consider $\alpha_1 = \alpha_2 = \alpha$, while for
a one-sided confidence interval of nominal level $1 - \alpha$ we would
consider either $\alpha_1 = 0$ and $\alpha_2 = \alpha$, or $\alpha_1 =
\alpha$ and $\alpha_2 = 0$. In many cases, it is possible to
replace the $\liminf_{n \rightarrow\infty}$ and $\geq$ in (\ref
{eqtwosided}) with $\lim_{n \rightarrow\infty}$ and $=$,
respectively. These results differ from those usually stated in the
literature in that they require the convergence to hold uniformly over
$\mathbf P$ instead of just pointwise over $\mathbf P$. The importance
of this stronger notion of convergence when applying these results is
discussed further below.

As we will see, the result (\ref{eqtwosided}) may hold with $\alpha_1 =
0$ and $\alpha_2 = \alpha\in(0,1)$, but it may fail if $\alpha_2 = 0$
and $\alpha_1 = \alpha\in(0,1)$, or the other way round.
This phenomenon arises when it is not possible to estimate $J_n(x,P)$
uniformly well with respect to a suitable metric, but, in a sense to be
made precise by our results, it is possible to estimate it sufficiently
well to ensure that (\ref{eqtwosided}) still holds for certain
choices of $\alpha_1$ and $\alpha_2$. Note that metrics compatible
with the weak topology are not sufficient for our purposes. In
particular, closeness of distributions with respect to such a metric
does not ensure closeness of quantiles. See Remark~\ref{remarklevy}
for further discussion of this point. In fact, closeness of
distributions with respect to even stronger metrics, such as the
Kolmogorov metric, does not ensure closeness of quantiles either. For
this reason, our results rely heavily on Lemma~\ref{lemmaquant} which
relates closeness of distributions with respect to a suitable metric
and coverage statements.

In contrast, the usual arguments for the pointwise asymptotic validity
of subsampling and the bootstrap rely on showing for each $P \in
\mathbf P$ that $\hat c_n ( 1- \alpha)$ tends in probability under $P$
to the $1- \alpha$ quantile of the limiting distribution of $R_n$
under~$P$. Because our results are uniform in $P \in{\mathbf P}$, we must
consider the behavior of $R_n$ and $\hat c_n ( 1- \alpha)$ under
arbitrary sequences $\{ P_n \in\mathbf P\dvtx  n \geq1\}$, under which
the quantile estimators need not even settle down. Thus, the results
are not trivial extensions of the usual pointwise asymptotic arguments.

The construction of $\hat c_n(\alpha)$ satisfying (\ref{eqtwosided})
is useful for constructing confidence regions that behave well
uniformly over $\mathbf P$. More precisely, our results provide
conditions under which subsampling and the bootstrap can be used to
construct confidence regions $C_n = C_n(X^{(n)})$ of level $1-\alpha$
for a parameter $\theta(P)$ that are uniformly consistent in level in
the sense that
%
\begin{equation}
\label{equnifcr} \liminf_{n \rightarrow\infty} \inf_{P \in\mathbf P}
P\bigl\{\theta(P)
\in C_n\bigr\} \geq1 -\alpha.
\end{equation}
Our results are also useful for constructing tests $\phi_n = \phi
_n(X^{(n)})$ of level $\alpha$ for a null hypothesis $P \in\mathbf
P_0 \subseteq\mathbf P$ against the alternative $P \in\mathbf P_1 =
\mathbf P \setminus\mathbf P_0$ that are uniformly consistent in level
in the sense that
%
\begin{equation}
\label{equniftest} \limsup_{n \rightarrow\infty} \sup_{P \in\mathbf
P_0} E_P[
\phi_n] \leq\alpha.
\end{equation}
In some cases, it is possible to replace the $\liminf_{n \rightarrow
\infty}$ and $\geq$ in (\ref{equnifcr}) or the $\limsup_{n
\rightarrow\infty}$ and $\leq$ in (\ref{equniftest}) with $\lim_{n
\rightarrow\infty}$ and $=$, respectively.

Confidence regions satisfying (\ref{equnifcr}) are desirable because
they ensure that for every $\varepsilon> 0$ there is an $N$ such that for
$n > N$ we have that $P\{\theta(P) \in C_n\}$ is no less than $1 -
\alpha- \varepsilon$ for all $P \in\mathbf P$. In contrast, confidence
regions that are only pointwise consistent in level in the sense that
\[
\liminf_{n \rightarrow\infty} P\bigl\{\theta(P) \in C_n\bigr\} \geq1 -
\alpha
\]
for each fixed $P \in\mathbf P$
have the feature that there exists some $\varepsilon> 0$ and $\{P_n \in
\mathbf P\dvtx  n \geq1 \}$ such that $P_n\{\theta(P_n) \in C_n\}$ is
less than $1 - \alpha- \varepsilon$ infinitely often. Likewise, tests
satisfying (\ref{equniftest}) are desirable for analogous reasons.
For this reason, inferences based on confidence regions or tests that
fail to satisfy (\ref{equnifcr}) or (\ref{equniftest}) may be very
misleading in finite samples. Of course, as pointed out by \citet
{bahadursavage1956}, there may be no nontrivial confidence region or
test satisfying (\ref{equnifcr}) or (\ref{equniftest}) when
$\mathbf P$ is sufficiently rich. For this reason, we will have to
restrict $\mathbf P$ appropriately in our examples. In the case of
confidence regions for or tests about the mean, for instance, we will
have to impose a very weak uniform integrability condition. See also
\citet{kabaila1995effect}, \citet{potscher2002lower},
Leeb and P{\"o}tscher
(\citeyear{leeb2006can,leeb2006performance}),
\citet{potscher2009confidence} for related results in more complicated
settings, including post-model selection, shrinkage-estimators and
ill-posed problems.

Some of our results on subsampling are closely related to results in
\citet{andrewsguggenberger2010}, which were developed independently
and at about the same time as our results. See the discussion on page
431 of \citet{andrewsguggenberger2010}. Our results show that the
question of whether subsampling can be used to construct estimators
$\hat c_n(\alpha)$ satisfying (\ref{eqtwosided}) reduces to a
single, succinct requirement on the asymptotic relationship between the
distribution of $J_n(x,P)$ and $J_b(x,P)$, where $b$ is the subsample
size, whereas the results of \citet{andrewsguggenberger2010} require
the verification of a larger number of conditions. Moreover, we also
provide a converse, showing this requirement on the asymptotic
relationship between the distribution of $J_n(x,P)$ and $J_b(x,P)$ is
also necessary
in the sense that, if the requirement fails, then for some nominal
coverage level, the uniform
coverage statements fail. Thus our results are stated under essentially
the weakest possible conditions, yet are verifiable in a large class of
examples. On the other hand, the results of \citet
{andrewsguggenberger2010} further provide a means of calculating the
limiting value of $\inf_{P \in\mathbf P} P\{\hat c_n(\alpha_1) \leq
R_n \leq\hat c_n(1 - \alpha_2)\}$ in the case where it may not
satisfy (\ref{eqtwosided}). To the best of our knowledge, our results
on the bootstrap are the first to be stated at this level of
generality. An important antecedent is \citet{romano1989unif}, who
studies the uniform asymptotic behavior of confidence regions for a
univariate cumulative distribution function. See also \citet
{mikusheva2007}, who analyzes the uniform asymptotic behavior of some
tests that arise in the context of an autoregressive model.

The remainder of the paper is organized as follows. In Section \ref
{sectiongeneral}, we present the conditions under which $\hat
c_n(\alpha)$ satisfying (\ref{eqtwosided}) may be constructed using
subsampling or the bootstrap. We then provide in Section \ref
{sectionappl} several applications of our general results. These
applications include the multivariate mean, testing moment
inequalities, multiple testing, the empirical process and
$U$-statistics. The discussion of $U$-statistics is especially
noteworthy because it highlights the fact that the assumptions required
for the uniform asymptotic validity of subsampling and the bootstrap
may differ. In particular, subsampling may be uniformly asymptotically
valid under conditions where, as noted by
\citet{bickelfreedman1981}, the bootstrap fails even to be
pointwise asymptotically valid. The application to multiple testing is
also noteworthy because, despite the enormous recent literature in this
area, our results appear to be the first that provide uniformly
asymptotically valid inference. Proofs of the main results (Theorems
\ref {theoremsubsample} and~\ref{theoremboot}) can be found in the
\hyperref[app]{Appendix}; proofs of all other results can be found in
\citet{romanoshaikh2012supp}, which contains supplementary
material. Many of the intermediate results may be of independent
interest, including uniform weak laws of large numbers for
$U$-statistics and $V$-statistics [Lemmas S.17.3 and S.17.4 in
\citet{romanoshaikh2012supp}, resp.] as well as the aforementioned
Lemma~\ref{lemmaquant}.

\section{General results} \label{sectiongeneral}

\subsection{Subsampling}

Let $X^{(n)} = (X_1,\ldots, X_n)$ be an i.i.d. sequence of random
variables with distribution $P \in\mathbf{P}$. Denote by $J_n(x,P)$
the distribution of a real-valued root $R_n = R_n(X^{(n)},P)$ under
$P$. The goal is to construct procedures which are valid uniformly in
$P$. In order to describe the subsampling approach to approximate
$J_n(x,P)$, let $b = b_n < n$ be a sequence of positive integers
tending to infinity, but satisfying $b/n \rightarrow0$, and define
$N_n = {n \choose b}$. For $i = 1,\ldots, N_n$, denote by
$X^{n,(b),i}$ the $i$th subset of data of size $b$. Below, we present
results for two subsampling-based estimators of $J_n(x,P)$. We first
consider the estimator given by
%
\begin{equation}
\label{equationsubdist} L_n(x,P) = \frac{1}{N_n} \sum
_{1 \leq i \leq N_n} I\bigl\{ R_b\bigl(X^{n,(b),i},P\bigr)
\leq x\bigr\}.
\end{equation}
More generally,
we will also consider feasible estimators $\hat L_n(x)$ in which $R_b$
is replaced by some estimator $\hat R_b$, that is,
%
\begin{equation}
\label{equationfeassubdist} \hat L_n(x) = \frac{1}{N_n} \sum
_{1 \leq i \leq N_n} I\bigl\{\hat R_b
\bigl(X^{n,(b),i}\bigr) \leq x\bigr\}.
\end{equation}
Typically, $\hat R_b(\cdot) = R_b(\cdot, \hat P_n)$, where $\hat P_n$
is the empirical distribution, but this is not assumed below. Even
though the estimator of $J_n(x,P)$ defined in (\ref{equationsubdist})
is infeasible because of its dependence on $P$, which is unknown, it is
useful both as an intermediate step toward establishing some results
for the feasible estimator of $J_n(x,P)$ and, as explained in Remarks
\ref{remarkinvert} and~\ref{remarktesting}, on its own in the
construction of some feasible tests and confidence regions.
%
\begin{theorem} \label{theoremsubsample}
Let $b = b_n < n$ be a sequence of positive integers tending to
infinity, but satisfying $b/n \rightarrow0$, and define $L_n(x,P)$ as
in (\ref{equationsubdist}). Then, the following statements are true:
\begin{longlist}[(iii)]
\item[(i)] If $\limsup_{n \rightarrow\infty} \sup_{P \in\mathbf
{P}} \sup_{x \in\mathbf{R}} \{J_b(x,P) - J_n(x,P)\} \leq0$, then
%
\begin{equation}
\label{eqtwosidedsub} \liminf_{n \rightarrow\infty} \inf_{P \in\mathbf
{P}} P\bigl\{
L_n^{-1}(\alpha_1,P) \leq R_n \leq
L_n^{-1}(1-\alpha_2,P)\bigr\} \geq1 -
\alpha_1 - \alpha_2
\end{equation}
holds for $\alpha_1 = 0$ and any $0 \leq\alpha_2 < 1$.
\item[(ii)] If $\limsup_{n \rightarrow\infty} \sup_{P \in\mathbf
{P}} \sup_{x \in\mathbf{R}} \{J_n(x,P) - J_b(x,P)\} \leq0$, then
(\ref{eqtwosidedsub}) holds for $\alpha_2 = 0$ and any $0 \leq
\alpha_1 < 1$.
\item[(iii)] If $\lim_{n \rightarrow\infty} \sup_{P \in\mathbf
{P}} \sup_{x \in\mathbf{R}} |J_b(x,P) - J_n(x,P)| = 0$, then (\ref
{eqtwosidedsub}) holds for any $\alpha_1 \geq0$ and $\alpha_2 \geq
0$ satisfying $0 \leq\alpha_1 + \alpha_2 < 1$.
\end{longlist}
\end{theorem}
%
\begin{remark} \label{remarkexactsub}
It is typically easy to deduce from the conclusions of Theorem~\ref
{theoremsubsample} stronger results in which the $\liminf_{n
\rightarrow\infty}$ and $\geq$ in (\ref{eqtwosidedsub}) are replaced by
$\lim_{n \rightarrow\infty}$ and $=$, respectively. For example, in
order to assert that (\ref{eqtwosidedsub}) holds with $\liminf_{n
\rightarrow\infty}$ and $\geq$ replaced by $\lim_{n \rightarrow\infty}$
and $=$, respectively, all that is required is that
\[
\lim_{n \rightarrow\infty} P\bigl\{L_n^{-1}(\alpha_1,P)
\leq R_n \leq L_n^{-1}(1 -
\alpha_2,P)\bigr\} = 1 - \alpha_1 - \alpha_2
\]
for some $P \in\mathbf P$. This can be verified using the usual
arguments for the pointwise asymptotic validity of subsampling. Indeed,
it suffices to show for some $P \in\mathbf P$ that $J_n(x,P)$ tends in
distribution to a limiting distribution $J(x,P)$ that is continuous at
the appropriate quantiles. See \citet{politisromanowolf1999} for
details.
\end{remark}
%
\begin{remark} \label{remarkinvert}
As mentioned earlier, $L_n(x,P)$ defined in (\ref{equationsubdist})
is infeasible because it still depends on $P$, which is unknown,
through $R_b(X^{n,(b),i},P)$. Even so, Theorem~\ref{theoremsubsample}
may be used without modification to construct feasible confidence
regions for a parameter of interest $\theta(P)$ provided that
$R_n(X^{(n)},P)$, and therefore $L_n(x,P)$, depends on $P$ only through
$\theta(P)$. If this is the case, then one may simply invert tests of
the null hypotheses $\theta(P) = \theta$ for all $\theta\in\Theta$
to construct a confidence region for $\theta(P)$. More concretely,
suppose $R_n(X^{(n)},P) = R_n(X^{(n)},\theta(P))$ and $L_n(x,P) =
L_n(x,\theta(P))$. Whenever we may apply part (i) of Theorem \ref
{theoremsubsample}, we have that
\[
C_n = \bigl\{\theta\in\Theta\dvtx  R_n\bigl(X^{(n)},
\theta\bigr) \leq L_n^{-1}(1 - \alpha, \theta)\bigr\}
\]
satisfies (\ref{equnifcr}). Similar conclusions follow from parts
(ii) and (iii) of Theorem~\ref{theoremsubsample}.
\end{remark}
%
\begin{remark} \label{remarktesting}
It is worth emphasizing that even though Theorem \ref
{theoremsubsample} is stated for roots, it is, of course, applicable
in the special case where $R_n(X^{(n)},P) = T_n(X^{(n)})$. This is
especially useful in the context of hypothesis testing. See Example
\ref{exmomineqsub} for one such instance.
\end{remark}

Next, we provide some results for feasible estimators of $J_n(x,P)$.
The first result,
Corollary~\ref{corrsubsample},
handles the case of the most basic root, while Theorem~\ref{theoremsubsample2}
applies to more general roots needed for many of our applications.

\begin{corollary} \label{corrsubsample}
Suppose $R_n = R_n(X^{(n)},P) = \tau_n(\hat\theta_n - \theta(P))$,
where $\{ \tau_n \in\mathbf R\dvtx  n \geq1\}$ is a sequence of
normalizing constants, $\theta(P)$ is a real-valued parameter of
interest and $\hat\theta_n = \hat\theta_n(X^{(n)})$ is an estimator
of $\theta(P)$. Let $b = b_n < n$ be a sequence of positive integers
tending to infinity, but satisfying $b/n \rightarrow0$, and define
\[
\hat L_n(x) = \frac{1}{N_n} \sum_{1 \leq i \leq N_n}
I\bigl\{\tau_b \bigl(\hat\theta_b\bigl(X^{n,(b),i}
\bigr) - \hat\theta_n\bigr) \leq x\bigr\}.
\]
Then statements \textup{(i)--(iii)} of Theorem~\ref{theoremsubsample} hold
when $L_n^{-1}(\cdot,P)$ is replaced by $\frac{\tau_n}{\tau_n +
\tau_b}\hat L_n^{-1}(\cdot)$.
\end{corollary}
%
\begin{theorem} \label{theoremsubsample2}
Let $b = b_n < n$ be a sequence of positive integers tending to
infinity, but satisfying $b/n \rightarrow0$. Define $L_n(x,P)$ as in
(\ref{equationsubdist}) and $\hat L_n(x)$ as in (\ref
{equationfeassubdist}). Suppose for all $\varepsilon> 0$ that
%
\begin{equation}
\label{eqfeasclose} \sup_{P \in\mathbf P} P \Bigl\{\sup_{x \in\mathbf
R} \bigl|\hat
L_n(x) - L_n(x,P)\bigr| > \varepsilon\Bigr\} \rightarrow0.
\end{equation}
Then, statements \textup{(i)--(iii)} of Theorem~\ref{theoremsubsample}
hold when $L_n^{-1}(\cdot,P)$ is replaced by $\hat L_n^{-1}(\cdot)$.
\end{theorem}

As a special case, Theorem~\ref{theoremsubsample2} can be applied to
Studentized roots.
%
\begin{corollary} \label{corrsubsample2}
Suppose
\[
R_n = R_n\bigl(X^{(n)},P\bigr) =
\frac{\tau_n(\hat\theta_n - \theta
(P))}{\hat\sigma_n},
\]
where $\{ \tau_n \in\mathbf R\dvtx  n \geq1\}$ is a sequence of
normalizing constants, $\theta(P)$ is a real-valued parameter of
interest, and $\hat\theta_n = \hat\theta_n(X^{(n)})$ is an
estimator of $\theta(P)$, and $\hat\sigma_n = \hat\sigma_n(X^{(n)}) \geq
0$ is an estimator of some parameter $\sigma(P) \geq
0$. Suppose further that:
\begin{longlist}[(ii)]
\item[(i)] The family of distributions $\{J_n(x,P)\dvtx  n \geq1, P \in
\mathbf P\}$ is tight, and any subsequential limiting distribution is
continuous.\vadjust{\goodbreak}
\item[(ii)] For any $\varepsilon> 0$,
\[
\sup_{P \in\mathbf P} P \biggl\{\biggl\llvert\frac{\hat\sigma_n}{\sigma
(P)} - 1 \biggr\rrvert
> \varepsilon\biggr\} \rightarrow0.
\]
\end{longlist}
Let $b = b_n < n$ be a sequence of positive integers tending to
infinity, but satisfying $b/n \rightarrow0$ and $\tau_b/\tau_n
\rightarrow0$. Define
\[
\hat L_n(x) = \frac{1}{N_n} \sum_{1 \leq i \leq N_n}
I \biggl\{\frac
{\tau_b (\hat\theta_b(X^{n,(b),i}) - \hat\theta_n)}{\hat\sigma
_b(X^{n,(b),i})} \leq x \biggr\}.
\]
Then statements \textup{(i)--(iii)} of Theorem~\ref{theoremsubsample} hold
when $L_n^{-1}(\cdot,P)$ is replaced by $\hat
L_n^{-1}(\cdot)$.\vspace*{-2pt}
\end{corollary}
%
\begin{remark}\label{remarkdec6}
One can take $\hat\sigma_n = \sigma(P)$ in Corollary~\ref{corrsubsample2}.
Since $\sigma(P)$ effectively cancels out from both sides of the
inequality in the event
$\{ R_n \le\hat L_n^{-1} ( 1- \alpha) \}$, such a root actually
leads to a computationally
feasible construction. However, Corollary~\ref{corrsubsample2} still
applies and shows
that we can obtain a positive result without the correction factor
$\tau_n / ( \tau_n + \tau_b )$ present in Corollary \ref
{corrsubsample}, provided the conditions of Corollary \ref
{corrsubsample2} hold. For example, if for some $\sigma(P)$, we have
that $\tau_n ( \hat\theta_n - \theta(P_n )) / \sigma(P_n )$ is
asymptotically standard normal under any sequence $\{P_n \in\mathbf P
\dvtx n \geq1\}$, then the conditions hold.\vspace*{-2pt}
\end{remark}
%
\begin{remark}
In Corollaries~\ref{corrsubsample} and~\ref{corrsubsample2}, it is
assumed that the rate of convergence $\tau_n$ is known. This
assumption may be relaxed using techniques described in \citet
{politisromanowolf1999}.\vspace*{-2pt}
\end{remark}

We conclude this section with a result that establishes a converse for
Theorems~\ref{theoremsubsample} and~\ref{theoremsubsample2}.\vspace*{-2pt}
%
\begin{theorem} \label{theoremconverse}
Let\vspace*{1pt} $b = b_n < n$ be a sequence of positive integers tending to
infinity, but satisfying $b/n \rightarrow0$ and define $L_n(x,P)$ as
in (\ref{equationsubdist}) and $\hat L_n(x)$ as in (\ref
{equationfeassubdist}). Then the following statements are true:
\begin{longlist}[(iii)]
\item[(i)] If $\limsup_{n \rightarrow\infty} \sup_{P \in\mathbf
{P}} \sup_{x \in\mathbf{R}} \{J_b(x,P) - J_n(x,P)\} > 0$, then (\ref
{eqtwosidedsub}) fails for $\alpha_1 = 0$ and some $0 \leq\alpha_2
< 1$.
\item[(ii)] If $\limsup_{n \rightarrow\infty} \sup_{P \in\mathbf
{P}} \sup_{x \in\mathbf{R}} \{J_n(x,P) - J_b(x,P)\} > 0$, then (\ref
{eqtwosidedsub}) fails for $\alpha_2 = 0$ and some $0 \leq\alpha_1
< 1$.
\item[(iii)] If $\liminf_{n \rightarrow\infty} \sup_{P \in\mathbf
{P}} \sup_{x \in\mathbf{R}} |J_b(x,P) - J_n(x,P)| > 0$, then (\ref
{eqtwosidedsub}) fails for some $\alpha_1 \geq0$ and $\alpha_2 \geq
0$ satisfying $0 \leq\alpha_1 + \alpha_2 < 1$.
\end{longlist}
If, in addition, (\ref{eqfeasclose}) holds for any $\varepsilon> 0$,
then statements \textup{(i)--(iii)} above hold when $ L_n^{-1}(\cdot,P)$ is
replaced by $\hat L_n^{-1}(\cdot)$.\vspace*{-2pt}
\end{theorem}

\subsection{Bootstrap}

As before, let $X^{(n)} = (X_1,\ldots, X_n)$ be an i.i.d. sequence of
random variables with distribution $P \in\mathbf{P}$. Denote by
$J_n(x,P)$ the distribution of a real-valued root $R_n =
R_n(X^{(n)},P)$ under $P$.
The goal remains to construct procedures which are valid\vadjust{\goodbreak} uniformly in $P$.
The bootstrap approach is to approximate $J_n ( \cdot, P )$ by
$J_n ( \cdot, \hat P_n )$ for some estimator $\hat P_n$ of $P$.
Typically, $\hat P_n$ is the empirical distribution, but this is not assumed
in Theorem~\ref{theoremboot} below.
Because $\hat P_n$ need not a priori even lie in ${\mathbf P}$, it is necessary
to introduce a family $\mathbf{P}'$ in which $\hat P_n$ lies (at least
with high probability). In order for the bootstrap to succeed,
we will require that $\rho( \hat P_n, P )$ be small for some
function (perhaps a metric) $\rho( \cdot, \cdot)$ defined
on $\mathbf{P}' \times\mathbf{P}$. For any given problem in
which the theorem is applied, $\mathbf{P}$, $\mathbf{P}'$ and $\rho$
must be specified.
%
\begin{theorem} \label{theoremboot}
Let $\rho(\cdot,\cdot)$ be a function on $\mathbf{P}' \times
\mathbf{P}$, and let $\hat P_n$ be a (random) sequence of
distributions. Then, the following are true:
\begin{longlist}[(iii)]
\item[(i)] Suppose $\limsup_{n \rightarrow\infty} \sup_{x \in
\mathbf{R}} \{J_n(x,Q_n) - J_n(x,P_n)\} \leq0$ for any sequences $\{
Q_n \in\mathbf{P}'\dvtx  n \geq1\}$ and $\{P_n \in\mathbf{P}\dvtx  n \geq
1\}$ satisfying $\rho(Q_n,P_n) \rightarrow0$. If
%
\begin{equation}
\label{eqihatethishere} \rho(\hat P_n, P_n)
\stackrel{P_n} {\rightarrow} 0 \quad\mbox{and}\quad P_n\bigl\{
\hat P_n \in\mathbf P'\bigr\} \rightarrow1
\end{equation}
for any sequence $\{P_n \in\mathbf P\dvtx  n \geq1\}$, then
%
\begin{equation}
\label{eqtwosidedboot} \liminf_{n \rightarrow\infty} \inf_{P \in\mathbf
{P}} P\bigl\{
J_n^{-1}(\alpha_1, \hat P_n) \leq
R_n \leq J_n^{-1}(1-\alpha_2, \hat
P_n)\bigr\} \geq1 - \alpha_1 - \alpha_2
\end{equation}
holds for $\alpha_1 = 0$ and any $0 \leq\alpha_2 < 1$.
\item[(ii)] Suppose $\limsup_{n \rightarrow\infty} \sup_{x \in
\mathbf{R}} \{J_n(x,P_n) - J_n(x,Q_n)\} \leq0$ for any sequences $\{
Q_n \in\mathbf{P}'\dvtx  n \geq1\}$ and $\{P_n \in\mathbf{P}\dvtx  n \geq 1\}$
satisfying $\rho(Q_n,P_n) \rightarrow0$. If (\ref{eqihatethishere})
holds for any sequence $\{P_n \in\mathbf P\dvtx  n \geq1\}$, then
(\ref{eqtwosidedboot}) holds for $\alpha_2 = 0$ and any $0 \leq\alpha_1
< 1$.
\item[(iii)] Suppose $\lim_{n \rightarrow\infty} \sup_{x \in
\mathbf{R}} |J_n(x,Q_n) - J_n(x,P_n)| = 0$ for any sequences $\{Q_n
\in\mathbf{P}'\dvtx  n \geq1\}$ and $\{P_n \in\mathbf{P}\dvtx  n \geq1\}$
satisfying $\rho(Q_n,P_n) \rightarrow0$. If (\ref{eqihatethishere})
holds for any sequence $\{P_n \in\mathbf P\dvtx  n \geq1\}$, then (\ref
{eqtwosidedboot}) holds for any $\alpha_1 \geq0$ and $\alpha_2 \geq
0$ satisfying $0 \leq\alpha_1 + \alpha_2 < 1$.
\end{longlist}
\end{theorem}
%
\begin{remark} \label{remarkexactboot}
It is typically easy to deduce from the conclusions of
Theorem~\ref{theoremboot} stronger results in which the $\liminf_{n
\rightarrow \infty}$ and $\geq$ in (\ref{eqtwosidedboot}) are replaced
by $\lim_{n \rightarrow\infty}$ and $=$, respectively. For example, in
order to assert that (\ref{eqtwosidedboot}) holds with $\liminf_{n
\rightarrow\infty}$ and $\geq$ replaced by $\lim_{n \rightarrow
\infty}$ and $=$, respectively, all that is required is that
\[
\lim_{n \rightarrow\infty} P\bigl\{J_n^{-1}(\alpha_1,
\hat P_n) \leq R_n \leq J_n^{-1}(1
- \alpha_2,\hat P_n)\bigr\} = 1 - \alpha_1 -
\alpha_2
\]
for some $P \in\mathbf P$. This can be verified using the usual
arguments for the pointwise asymptotic validity of the bootstrap. See
\citet{politisromanowolf1999} for details.
\end{remark}
%
\begin{remark} \label{remarklevy}
In some cases, it is possible to construct estimators $\hat J_n(x)$ of
$J_n(x,P)$ that are uniformly consistent over a large class of
distributions $\mathbf P$ in the sense that for any $\varepsilon> 0$
%
\begin{equation}
\label{eqlevyconv} \sup_{P \in\mathbf P} P\bigl\{\rho\bigl(\hat J_n(
\cdot), J_n(\cdot,P)\bigr) > \varepsilon\bigr\} \rightarrow0,
\end{equation}
where $\rho$ is the Levy metric or some other metric compatible with
the weak topology. Yet a result such as (\ref{eqlevyconv}) is not
strong enough to yield uniform coverage statements such as those in
Theorems~\ref{theoremsubsample} and~\ref{theoremboot}. In
other words, such conclusions do not follow from uniform approximations
of the distribution of interest if the quality of the approximation is
measured in terms of metrics metrizing weak convergence. To see this,
consider the following simple example.
%
\begin{example}
Let $X^{(n)} = (X_1,\ldots, X_n)$ be an i.i.d. sequence of random
variables with distribution $P_\theta= \operatorname
{Bernoulli}(\theta)$.
Denote by $J_n(x,P_\theta)$ the distribution of the root $R_n = \sqrt n
(\hat\theta_n - \theta)$ under $P_\theta$, where $\hat\theta_n =
\bar
X_n$. Let $\hat P_n$ be the empirical distribution of $X^{(n)}$ or,
equivalently, $P_{\hat\theta_n}$. Lemma S.1.1 in
\citet{romanoshaikh2012supp} implies for any $\varepsilon> 0$ that
%
\begin{equation}
\label{eqlevy} \sup_{0 \le\theta\le1} P_\theta\bigl\{\rho
\bigl(J_n(\cdot,\hat P_n), J_n(
\cdot,P_\theta)\bigr) > \varepsilon\bigr\} \rightarrow0,
\end{equation}
whenever $\rho$ is a metric compatible with the weak topology.
Nevertheless, it follows from the argument on page 78 of \citet
{romano1989unif} that the coverage statements in Theorem \ref
{theoremboot} fail to hold provided that both $\alpha_1$ and $\alpha_2$
do not equal zero. Indeed, consider part (i) of Theorem \ref
{theoremboot}. Suppose $\alpha_1 = 0$ and $0 < \alpha_2 < 1$. For a
given $n$ and $\delta> 0$, let $\theta_n = (1 - \delta)^{
{1}/{n}}$. Under $P_{\theta_n}$, the event $X_1 = \cdots= X_n = 1$ has
probability $1 - \delta$. Moreover, whenever such an event occurs,
$R_n > J_n^{-1}(1 - \alpha_2, \hat P_n) = 0$. Therefore, $P_{\theta
_n}\{J_n^{-1}(\alpha_1,\hat P_n) \leq R_n \leq J_n^{-1}(1 - \alpha
_2,\hat P_n)\} \leq\delta$. Since the choice of $\delta$ was
arbitrary, it follows that
\[
\liminf_{n \rightarrow\infty} \inf_{0 \leq\theta\leq1} P_\theta\bigl
\{J_n^{-1}(\alpha_1,\hat P_n) \leq
R_n \leq J_n^{-1}(1 - \alpha_2,
\hat P_n)\bigr\} = 0.
\]
A similar argument establishes the result for parts (ii) and (iii) of
Theorem~\ref{theoremboot}.
\end{example}

On the other hand, when $\rho$ is the Kolmogorov metric,
(\ref{eqlevy}) holds when the supremum over $0 \le\theta\le1$ is
replaced with a supremum over $\delta< \theta< 1 - \delta$ for some
$\delta> 0$. Moreover, when $\theta$ is restricted to such an
interval, the coverage statements in Theorem~\ref{theoremboot} hold
as well.
\end{remark}

\section{Applications} \label{sectionappl}

Before proceeding, it is useful to introduce some notation that will be
used frequently throughout many of the examples below. For a
distribution $P$ on $\mathbf R^k$, denote by $\mu(P)$ the mean of $P$,
by $\Sigma(P)$ the covariance matrix of~$P$, and by $\Omega(P)$ the
correlation matrix of $P$. For $1 \leq j \leq k$, denote by $\mu_j(P)$
the $j$th component of $\mu(P)$ and by $\sigma^2_j(P)$ the $j$th
diagonal\vspace*{-1pt} element of $\Sigma(P)$. In all of our examples, $X^{(n)} =
(X_1,\ldots, X_n)$ will\vspace*{1pt} be an i.i.d. sequence of random variables
with distribution $P$ and $\hat P_n$ will denote the empirical
distribution of $X^{(n)}$. As usual, we will denote by $\bar X_n = \mu
(\hat P_n)$ the usual sample mean, by $\hat\Sigma_n = \Sigma(\hat
P_n)$ the usual sample covariance matrix and by $\hat\Omega_n =
\Omega(\hat P_n)$ the usual sample correlation matrix.\vadjust{\goodbreak} For $1 \leq j
\leq k$, denote by $\bar X_{j,n}$ the $j$th component of $\bar X_n$ and
by $S^2_{j,n}$ the $j$th diagonal element of $\hat\Sigma_n$. Finally,
we say that a family of distributions $\mathbf Q$ on the real line
satisfies the standardized uniform integrability condition if
%
\begin{equation}
\label{equi} \lim_{\lambda\rightarrow\infty} \sup_{Q \in\mathbf Q} E_Q
\biggl[
\biggl(\frac{Y - \mu(Q)}{\sigma(Q)} \biggr)^2 I \biggl\{ \biggl\llvert
\frac{Y - \mu(Q)}{\sigma(Q)} \biggr\rrvert> \lambda\biggr\} \biggr] = 0.
\end{equation}
In the preceding expression, $Y$ denotes a random variable with
distribution $Q$. The use of the term standardized to describe (\ref
{equi}) reflects that fact that the variable $Y$ is centered around
its mean and normalized by its standard deviation.

\subsection{Subsampling}
%
\begin{example}[(Multivariate nonparametric mean)] \label{examplesubmean}
Let $X^{(n)} = (X_1,\ldots,\break X_n)$ be an i.i.d. sequence of random
variables with distribution $P \in\mathbf P$ on $\mathbf R^k$. Suppose
one wishes to construct a rectangular confidence region for $\mu(P)$.
For this purpose, a natural choice of root is
%
\begin{equation}
\label{eqmaxroot} R_n\bigl(X^{(n)},P\bigr) =
\max_{1 \leq j \leq k} \frac{\sqrt{n}(\bar X_{j,n}
- \mu_j(P))}{S_{j,n}}.
\end{equation}
In this setup, we have the following theorem:
%
\begin{theorem} \label{theoremsubmean}
Denote by $\mathbf P_j$ the set of distributions formed from the $j$th
marginal distributions of the distributions in $\mathbf P$. Suppose
$\mathbf P$ is such that (\ref{equi}) is satisfied with $\mathbf Q =
\mathbf P_j$ for all $1 \leq j \leq k$. Let $J_n(x,P)$ be the
distribution of the root (\ref{eqmaxroot}). Let $b = b_n < n$ be a
sequence of positive integers tending to infinity, but satisfying $b/n
\rightarrow0$ and define $L_n(x,P)$ by (\ref{equationsubdist}). Then
%
\begin{eqnarray}
\label{eqsubmean}\qquad
&&
\lim_{n \rightarrow\infty} \inf_{P \in\mathbf{P}} P
\biggl\{
L_n^{-1}(\alpha_1,P) \leq\max_{1 \leq j \leq k}
\frac{\sqrt{n}(\bar X_{j,n} - \mu_j(P))}{S_{j,n}} \leq
L_n^{-1}(1-\alpha_2,P)
\biggr\}\nonumber\\[-8pt]\\[-8pt]
&&\qquad = 1 - \alpha_1 - \alpha_2\nonumber
\end{eqnarray}
for any $\alpha_1 \geq0$ and $\alpha_2 \geq0$ such that $0 \leq
\alpha_1 + \alpha_2 < 1$. Furthermore, (\ref{eqsubmean}) remains
true if $L_n^{-1}(\cdot,P)$ is replaced by $\hat L_n^{-1}(\cdot)$,
where $\hat L_n(x)$ is defined by (\ref{equationfeassubdist}) with
$\hat R_b(X^{n,(b),i}) = R_b(X^{n,(b),i}, \hat P_n)$.
\end{theorem}

Under suitable restrictions, Theorem~\ref{theoremsubmean}
generalizes to the case where the root is given by
%
\begin{equation}
\label{eqgeneralf} R_n\bigl(X^{(n)},P\bigr) = f
\bigl(Z_n(P),\hat\Omega_n\bigr),
\end{equation}
where $f$ is a continuous, real-valued function and
%
\begin{equation}
\label{eqznp} Z_n(P) = \biggl(\frac{\sqrt{n}(\bar X_{1,n} - \mu
_1(P))}{S_{1,n}},\ldots,
\frac{\sqrt{n}(\bar X_{k,n} - \mu_k(P))}{S_{k,n}} \biggr)'.
\end{equation}
In particular, we have the following theorem:
%
\begin{theorem} \label{theoremsubgenmean}
Let $\mathbf P$ be defined as in Theorem~\ref{theoremsubmean}. Let
$J_n(x,P)$ be the distribution of root (\ref{eqgeneralf}), where $f$
is continuous.
\begin{longlist}[(ii)]
\item[(i)] Suppose further that for all $x \in\mathbf R$ that
%
\begin{eqnarray}
\label{eqgensubmeanconv1}
P_n\bigl\{f\bigl(Z_n(P_n),\Omega(\hat
P_n)\bigr) \leq x\bigr\} &\rightarrow& P\bigl\{ f(Z,\Omega) \leq x
\bigr\},
\\
\label{eqgensubmeanconv2}
P_n\bigl\{f\bigl(Z_n(P_n),\Omega(\hat
P_n)\bigr) < x\bigr\} &\rightarrow& P\bigl\{f(Z,\Omega) < x\bigr\}
\end{eqnarray}
for any sequence $\{P_n \in\mathbf P\dvtx  n \geq1\}$ such that $Z_n(P_n)
\stackrel{d}{\rightarrow} Z$ under $P_n$ and $\Omega(\hat P_n)
\stackrel{P_n}{\rightarrow} \Omega$, where $Z \sim N(0,\Omega)$. Then
%
\begin{eqnarray}
\label{eqgensubmeancover}
&&
\liminf_{n \rightarrow\infty} \inf_{P \in
\mathbf P} P\bigl\{
L_n^{-1}(\alpha_1,P) \leq f
\bigl(Z_n(P),\hat\Omega_n\bigr) \leq
L_n^{-1}(1 - \alpha_2,P)\bigr\} \nonumber\\[-8pt]\\[-8pt]
&&\qquad\geq1 -
\alpha_1 - \alpha_2\nonumber
\end{eqnarray}
for any $\alpha_1 \geq0$ and $\alpha_2 \geq0$ such that $0 \leq
\alpha_1 + \alpha_2 < 1$.
\item[(ii)] Suppose further that if $Z \sim N(0,\Omega)$ for some
$\Omega$ satisfying $\Omega_{j,j} = 1$ for all $1 \leq j \leq k$,
then $f(Z,\Omega)$ is continuously distributed. Then, (\ref
{eqgensubmeancover}) remains true if $L_n^{-1}(\cdot,P)$ is replaced
by $\hat L_n^{-1}(\cdot)$, where $\hat L_n(x)$ is defined by (\ref
{equationfeassubdist}) with $\hat R_b(X^{n,(b),i}) = R_b(X^{n,(b),i},
\hat P_n)$. Moreover, the $\liminf_{n \rightarrow\infty}$ and $\geq
$ may be replaced by $\lim_{n \rightarrow\infty}$ and $=$, respectively.
\end{longlist}
\end{theorem}

In order to verify (\ref{eqgensubmeanconv1}) and (\ref
{eqgensubmeanconv2}) in Theorem~\ref{theoremsubgenmean}, it suffices
to assume that $f(Z,\Omega)$ is continuously distributed. Under the
assumptions of the theorem, however, $f(Z,\Omega)$ need not be
continuously distributed. In this case, (\ref{eqgensubmeanconv1}) and
(\ref{eqgensubmeanconv2}) hold immediately for any $x$ at which $P\{
(Z,\Omega) \leq x\}$ is continuous, but require a further argument for
$x$ at which $P\{(Z,\Omega) \leq x\}$ is discontinuous. See, for
example, the proof of Theorem~\ref{theorembootmomineq}, which relies
on Theorem~\ref{theorembootgenmean}, where the same requirement
appears.
\end{example}
%
\begin{example}[(Constrained univariate nonparametric mean)]
\label{remarkonesided}
\citet{andrews2000} considers the following example. Let $X^{(n)} =
(X_1,\ldots, X_n)$ be an i.i.d. sequence of random variables with
distribution $P \in\mathbf P$ on $\mathbf R$. Suppose it is known that
$\mu(P) \geq0$ for all $P \in\mathbf P$ and one wishes to construct
a confidence interval for $\mu(P)$. A natural choice of root in this
case is
\[
R_n = R_n\bigl(X^{(n)},P\bigr) = \sqrt n\bigl(
\max\{\bar X_n,0\} - \mu(P)\bigr).
\]
This root differs from the one considered in Theorem \ref
{theoremsubmean} and the ones discussed in Theorem \ref
{theoremsubgenmean} in the sense that under weak assumptions on
$\mathbf P$,
%
\begin{equation}
\label{eqonesidedunif} \limsup_{n \rightarrow\infty} \sup_{P \in\mathbf{P}}
\sup_{x \in
\mathbf{R}} \bigl\{J_b(x,P) - J_n(x,P)\bigr\}
\leq0
\end{equation}
holds, but
%
\begin{equation}
\label{eqothersideunif} \limsup_{n \rightarrow\infty} \sup_{P \in\mathbf{P}}
\sup_{x \in
\mathbf{R}} \bigl\{J_n(x,P) - J_b(x,P)\bigr\}
\leq0
\end{equation}
fails to hold. To see this, suppose (\ref{equi}) holds with $\mathbf
Q = \mathbf P$. Note that
\begin{eqnarray*}
J_b(x,P) &=& P\bigl\{\max\bigl\{Z_b(P),-\sqrt b \mu(P)
\bigr\} \leq x\bigr\},
\\
J_n(x,P) &=& P\bigl\{\max\bigl\{Z_n(P),-\sqrt n \mu(P)
\bigr\} \leq x\bigr\},
\end{eqnarray*}
where $Z_b(P) = \sqrt b (\bar X_b - \mu(P))$ and $Z_n(P) = \sqrt n
(\bar X_n - \mu(P))$. Since $\sqrt b \mu(P) \leq\sqrt n \mu(P)$ for
any $P \in\mathbf P$, $J_b(x,P) - J_n(x,P)$ is bounded from above by
\[
P\bigl\{\max\bigl\{Z_b(P),-\sqrt n \mu(P)\bigr\} \leq x\bigr\} -
J_n(x,P).
\]
It now follows from the uniform central limit theorem established by
Lem\-ma~3.3.1 of \citet{romanoshaikh2008} and Theorem 2.11 of \citet
{bhattacharyarao1976} that (\ref{eqonesidedunif}) holds. It
therefore follows from Theorem~\ref{theoremsubsample} that (\ref
{eqtwosidedsub}) holds with $\alpha_1 = 0$ and any $0 \leq\alpha_2
< 1$. To see that (\ref{eqothersideunif}) fails, suppose further that
$\{Q_n\dvtx  n \geq1\} \subseteq\mathbf P$, where $Q_n = N(h/\sqrt n, 1)$
for some $h > 0$. For $Z \sim N(0,1)$,
\begin{eqnarray*}
J_n ( x, Q_n ) &=& P \bigl\{ \max(Z, -h ) \le x
\bigr\},
\\
J_b ( x, Q_n ) &=& P \bigl\{ \max(Z, -h \sqrt{b} /
\sqrt{n} ) \le x \bigr\}.
\end{eqnarray*}
The left-hand side of (\ref{eqothersideunif}) is therefore greater
than or equal to
\[
\limsup_{n \to\infty} \bigl( P \bigl\{ \max(Z, -h ) \le x \bigr\} - P
\bigl\{
\max(Z, -h \sqrt{b} / \sqrt{n} ) \le x \bigr\} \bigr)
\]
for any $x$. In particular, if $-h < x < 0$, then the second term is
zero for large enough $n$, and so the limiting value is $ P \{ Z \le x
\} = \Phi( x ) > 0$. It therefore follows from Theorem \ref
{theoremconverse} that (\ref{eqtwosidedsub}) fails for $\alpha_2 =
0$ and some $0 \leq\alpha_1 < 1$. On the other hand, (\ref
{eqtwosidedsub}) holds with $\alpha_2 = 0$ and any $0.5 < \alpha_1 <
1$. To see this, consider any sequence $\{P_n \in\mathbf P\dvtx  n\geq1\}
$ and the event $\{ L_n^{-1} ( \alpha_1, P_n ) \le R_n \}$. For the
root in this example, this event is scale invariant. So, in calculating
the probability of this event, we may without loss of generality assume
$\sigma^2 (P_n ) = 1$. Since $\mu( P_n ) \ge0$, we have for any $x
\ge0$ that
\[
J_n ( x, P_n ) = P \bigl\{ \max\bigl\{
Z_n ( P_n ), - \sqrt{n} \mu( P_n ) \bigr
\} \le x \bigr\} = P \bigl\{ Z_n ( P_n ) \le x \bigr\}
\to\Phi( x)
\]
and similarly for $J_b (x, P_n )$.
Using the usual subsampling arguments, it is thus possible to show for
$0.5 < \alpha_1 < 1$ that
\[
L_n^{-1} ( \alpha_1, P_n )
\stackrel{P_n} {\rightarrow} \Phi^{-1} (
\alpha_1 ).
\]
The desired conclusion therefore follows from Slutsky's theorem.
Arguing as the the proof of Corollary~\ref{corrsubsample2} and Remark
\ref{remarkdec6},
it can be shown that the same results hold
when $L_n^{-1}(\cdot,P)$ is replaced by $\hat L_n^{-1}(\cdot)$, where
$\hat L_n(x)$ is defined as $L_n(x,P)$ is defined but with $\mu(P)$
replaced by $\bar X_n$.
\end{example}
%
\begin{example}[(Moment inequalities)] \label{exmomineqsub}
The generality of Theorem~\ref{theoremsubsample} illustrated in Example
\ref{remarkonesided} is also useful when testing multisided hypotheses
about the mean. To see this, let $X^{(n)} = (X_1,\ldots, X_n)$ be an
i.i.d. sequence of random variables with distribution $P \in\mathbf P$
on $\mathbf R^k$. Define $\mathbf P_0 = \{P \in\mathbf P\dvtx  \mu(P)
\leq0\}$ and $\mathbf P_1 = \mathbf P \setminus\mathbf P_0$. Consider
testing the null hypothesis that $P \in\mathbf P_0$ versus the
alternative hypothesis that $P \in\mathbf P_1$ at level
$\alpha\in(0,1)$. Such hypothesis testing problems have recently
received considerable attention in the ``moment inequality'' literature
in econometrics. See, for example, \citet{andrewssoares2010},
\citet {andrewsguggenberger2010}, \citet{andrewsbarwick2012},
\citet {bugni2010}, \citet{canay2010} and Romano and Shaikh
(\citeyear{romanoshaikh2008,romanoshaikh2010}). Theorem
\ref{theoremsubsample} may be used to construct tests that are
uniformly consistent in level in the sense that (\ref{equniftest})
holds under weak assumptions on $\mathbf P$. Formally, we have the
following theorem:
%
\begin{theorem} \label{theoremsubmomineq}
Let $\mathbf P$ be defined as in Theorem~\ref{theoremsubmean}. Let
$J_n(x,P)$ be the distribution of
\[
T_n\bigl(X^{(n)}\bigr) = \max_{1 \leq j \leq k}
\frac{\sqrt n \bar X_{j,n}}{S_{j,n}}.
\]
Let $b = b_n < n$ be a sequence of positive integers tending to
infinity, but satisfying $b/n \rightarrow0$ and define $L_n(x)$ by the
right-hand side of (\ref{equationsubdist}) with $R_n(X^{(n)},P) =
T_n(X^{(n)})$. Then, the test defined by
\[
\phi_n\bigl(X^{(n)}\bigr) = I\bigl\{ T_n
\bigl(X^{(n)}\bigr) > L_n^{-1}(1 - \alpha)\bigr\}
\]
satisfies (\ref{equniftest}) for any $0 < \alpha< 1$.
\end{theorem}

The argument used to establish Theorem \ref
{theoremsubmomineq} is essentially the same as the one presented in
\citet{romanoshaikh2008} for
\[
T_n\bigl(X^{(n)}\bigr) = \sum_{1 \leq j \leq k}
\max\{\sqrt n \bar X_{j,n},0\}^2,
\]
though Lemma S.6.1 in \citet{romanoshaikh2012supp} is needed for
establishing (\ref{eqonesidedunif}) here because of Studentization.
Related results are obtained by \citet{andrewsguggenberger2009}.
\end{example}
%
\begin{example}[(Multiple testing)] \label{exmultsub}
We now illustrate the use of Theorem~\ref{theoremsubsample} to
construct tests of multiple hypotheses that behave well uniformly over
a large class of distributions.
Let $X^{(n)} = (X_1,\ldots, X_n)$ be an i.i.d. sequence of random
variables with distribution $P \in\mathbf P$ on $\mathbf R^k$, and
consider testing the family of null hypotheses
%
\begin{equation}
\label{eqfamily} H_j\dvtx  \mu_j(P) \leq0 \qquad\mbox{for } 1
\leq j \leq k
\end{equation}
versus the alternative hypotheses
%
\begin{equation}
\label{eqaltfamily} H_j'\dvtx  \mu_j(P) > 0
\qquad\mbox{for } 1 \leq j \leq k\vadjust{\goodbreak}
\end{equation}
in a way that controls the familywise error rate at level $0 < \alpha<
1$ in the sense that
%
\begin{equation}
\label{equniffwer} \limsup_{n \rightarrow\infty} \sup_{P \in\mathbf P}
\mathrm{FWER}_P \leq\alpha,
\end{equation}
where
\[
\mathrm{FWER}_P = P\bigl\{\mbox{reject some } H_j \mbox{ with }
\mu_j(P) \leq0\bigr\}.
\]
For $K \subseteq\{1,\ldots, k\}$, define $L_n(x,K)$ according to the
right-hand side of (\ref{equationsubdist}) with
\[
R_n\bigl(X^{(n)},P\bigr) = \max_{j \in K}
\frac{\sqrt n \bar X_{j,n}}{S_{j,n}},
\]
and consider the following stepwise multiple testing procedure:
%
\begin{algorithm} \label{algstepsub}
Step 1: Set $K_1 = \{1,\ldots, k\}$. If
\[
\max_{j \in K_1} \frac{\sqrt n \bar X_{j,n}}{S_{j,n}} \leq L_n^{-1}(1 -
\alpha,K_1),
\]
then stop. Otherwise, reject any $H_j$ with
\[
\frac{\sqrt n \bar X_{j,n}}{S_{j,n}} > L_n^{-1}(1 - \alpha,K_1)
\]
and continue to Step $2$ with
\[
K_{2} = \biggl\{j \in K_1\dvtx  \frac{\sqrt n \bar X_{j,n}}{S_{j,n}} \leq
L_n^{-1}(1 - \alpha,K_1) \biggr\}.
\]

$\vdots$\vspace*{10pt}

Step $s$: If
\[
\max_{j \in K_s} \frac{\sqrt n \bar X_{j,n}}{S_{j,n}} \leq L_n^{-1}(1 -
\alpha,K_s),
\]
then stop. Otherwise, reject any $H_j$ with
\[
\frac{\sqrt n \bar X_{j,n}}{S_{j,n}} > L_n^{-1}(1 - \alpha,K_s)
\]
and continue to Step $s+1$ with
\[
K_{s+1} = \biggl\{j \in K_s\dvtx  \frac{\sqrt n \bar X_{j,n}}{S_{j,n}} \leq
L_n^{-1}(1 - \alpha,K_s) \biggr\}.
\]

$\vdots$
\end{algorithm}

We have the following theorem:
%
\begin{theorem} \label{theoremsubmultiple}
Let $\mathbf P$ be defined as in Theorem~\ref{theoremsubmean}. Let $b
= b_n < n$ be a sequence of positive integers tending to infinity, but
satisfying $b/n \rightarrow0$. Then, Algorithm~\ref{algstepsub} satisfies
%
\begin{equation}
\label{eqfwercontrol} \limsup_{n \rightarrow\infty} \sup_{P \in\mathbf P}
\mathrm{FWER}_P \leq\alpha
\end{equation}
for any $0 < \alpha< 1$.
\end{theorem}

It is, of course, possible to extend the analysis in a straightforward
way to two-sided testing. See also \citet{romanoshaikh2010} for
related results about a multiple testing problem involving an infinite
number of null hypotheses.
\end{example}
%
\begin{example}[(Empirical process on $\mathbf R$)] \label{exsubemp}
Let $X^{(n)} = (X_1,\ldots, X_n)$ be an i.i.d. sequence of random
variables with distribution $P \in\mathbf P$ on $\mathbf R$. Suppose
one wishes to construct a confidence region for the cumulative
distribution function associated with $P$, that is, $P\{(-\infty,t]\}$.
For this purpose a natural choice of root is
%
\begin{equation}
\label{eqks} \sup_{t \in\mathbf R} \sqrt n \bigl|\hat P_n\bigl\{(-
\infty,t]\bigr\} - P\bigl\{(-\infty,t]\bigr\}\bigr|.
\end{equation}
In this setting, we have the following theorem:
%
\begin{theorem} \label{theoremsubemp}
Fix any $\varepsilon\in(0,1)$, and let
%
\begin{equation}
\label{eqnoatom} \mathbf P = \bigl\{P \mbox{ on } \mathbf R\dvtx  \varepsilon<
P\bigl
\{(- \infty,t]\bigr\} < 1 - \varepsilon\mbox{ for some } t \in\mathbf
R\bigr\}.
\end{equation}
Let $J_n(x,P)$ be the distribution of root (\ref{eqks}). Then
%
\begin{eqnarray}
\label{eqsubemp}
&&
\lim_{n \rightarrow\infty} \inf_{P \in\mathbf{P}} P
\Bigl\{
L_n^{-1}(\alpha_1, P) \leq\sup_{t \in\mathbf R}
\sqrt n \bigl|\hat P_n\bigl\{ (-\infty,t]\bigr\} - P\bigl\{(-\infty,t]
\bigr\}\bigr|\nonumber\\
&&\hspace*{182pt}\qquad \leq L_n^{-1}(1-\alpha_2, P) \Bigr\} \\
&&\qquad=
1 - \alpha_1 - \alpha_2\nonumber
\end{eqnarray}
for any $\alpha_1 \geq0$ and $\alpha_2 \geq0$ such that $0 \leq
\alpha_1 + \alpha_2 < 1$. Furthermore, (\ref{eqsubemp}) remains
true if $L_n^{-1}(\cdot,P)$ is replaced by $\hat L_n^{-1}(\cdot)$,
where $\hat L_n(x)$ is defined by (\ref{equationfeassubdist}) with
$\hat R_b(X^{n,(b),i}) = R_b(X^{n,(b),i}, \hat P_n)$.
\end{theorem}
\end{example}
%
\begin{example}[(One sample $U$-statistics)] \label{exsubustat}
Let $X^{(n)} = (X_1,\ldots, X_n)$ be an i.i.d. sequence of random
variables with distribution $P \in\mathbf P$ on $\mathbf R$. Suppose
one wishes to construct a confidence region for
%
\begin{equation}
\label{eqthetaustat} \theta(P) = \theta_h(P) = E_P
\bigl[h(X_1,\ldots,X_m)\bigr],
\end{equation}
where $h$ is a symmetric kernel of degree $m$. The usual estimator of
$\theta(P)$ in this case is given by the $U$-statistic
\[
\label{eqhatthetaustat} \hat\theta_{n} = \hat\theta_{n}
\bigl(X^{(n)}\bigr) = \frac{1}{{n \choose m}} \sum
_c h(X_{i_1},\ldots,X_{i_m}).
\]
Here, $\sum_c$ denotes summation over all ${n \choose m}$ subsets $\{
i_1,\ldots,i_m\}$ of $\{1,\ldots, n\}$. A~natural choice of root is
therefore given by
%
\begin{equation}
\label{equstatroot} R_n\bigl(X^{(n)},P\bigr) = \sqrt n
\bigl(\hat\theta_{n} - \theta(P)\bigr).
\end{equation}
In this setting, we have the following theorem:
%
\begin{theorem} \label{theoremsubustat}
Let
%
\begin{equation}
\label{eqgP} g (x,P) = g_h(x,P) = E_P \bigl[ h (x,
X_2,\ldots, X_m ) \bigr] - \theta(P)
\end{equation}
and
%
\begin{equation}
\label{eqsigmaP} \sigma^2_h(P) = m^2
\operatorname{Var}_P \bigl[ g (X_i, P ) \bigr].
\end{equation}
Suppose $\mathbf P$ satisfies the uniform integrability condition
%
\begin{equation}
\label{eqUUI} \lim_{\lambda\to\infty} \sup_{P \in{\mathbf P}} E_P \biggl[
\frac
{g^2 (X_i, P )}{\sigma_h^2 (P)} I \biggl\{ \biggl\llvert\frac
{g(X_i,P)}{\sigma_h (P)} \biggr\rrvert>
\lambda\biggr\} \biggr] = 0
\end{equation}
and
%
\begin{equation}
\label{eqerror} \sup_{P \in{\mathbf P }} \frac{\operatorname{Var}_P [ h (
X_1,\ldots, X_m ) ]}{\sigma^2 (P)} < \infty.
\end{equation}
Let $J_n(x,P)$ be the distribution of the root (\ref{equstatroot}).
Let $b = b_n < n$ be a sequence of positive integers tending to
infinity, but satisfying $b/n \rightarrow0$, and define $L_n(x,P)$ by
(\ref{equationsubdist}). Then
%
\begin{eqnarray}
\label{eqsubustat}
&&
\lim_{n \rightarrow\infty} \inf_{P \in\mathbf{P}} P
\bigl\{
L_n^{-1}(\alpha_1,P) \leq\sqrt n\bigl(\hat
\theta_{n} - \theta(P)\bigr) \leq L_n^{-1}(1-
\alpha_2,P) \bigr\}\nonumber\\[-8pt]\\[-8pt]
&&\qquad = 1 - \alpha_1 - \alpha_2\nonumber
\end{eqnarray}
for any $\alpha_1 \geq0$ and $\alpha_2 \geq0$ such that $0 \leq
\alpha_1 + \alpha_2 < 1$. Furthermore, (\ref{eqsubustat}) remains
true if $L_n^{-1}(\cdot,P)$ is replaced by $\hat L_n^{-1}(\cdot)$,
where $\hat L_n(x)$ is defined by (\ref{equationfeassubdist}) with
$\hat R_b(X^{n,(b),i}) = R_b(X^{n,(b),i}, \hat P_n)$.
\end{theorem}
\end{example}

\subsection{Bootstrap}
%
\begin{example}[(Multivariate nonparametric mean)] \label{examplebootmean}
Let $X^{(n)} = (X_1,\ldots,\break X_n)$ be an i.i.d. sequence of random
variables with distribution $P \in\mathbf P$ on $\mathbf R^k$. Suppose
one wishes to construct a rectangular confidence region for $\mu(P)$.
As described in Example~\ref{examplesubmean}, a natural choice of
root in this case is given by (\ref{eqmaxroot}). In this setting, we
have the following theorem, which is a bootstrap counterpart to Theorem
\ref{theoremsubmean}:
%
\begin{theorem} \label{theorembootmean}
Let $\mathbf P$ be defined as in Theorem~\ref{theoremsubmean}. Let
$J_n(x,P)$ be the distribution of the root (\ref{eqmaxroot}). Then
%
\begin{eqnarray}
\label{eqbootmeancover}
&&
\lim_{n \rightarrow\infty} \inf_{P \in\mathbf
{P}} P \biggl\{
J_n^{-1}(\alpha_1, \hat P_n) \leq
\max_{1 \leq j \leq k} \frac{\sqrt{n}(\bar X_{j,n} - \mu
_j(P))}{S_{j,n}} \leq J_n^{-1}(1-
\alpha_2, \hat P_n) \biggr\} \nonumber\hspace*{-35pt}\\[-8pt]\\[-8pt]
&&\qquad= 1 - \alpha_1 -
\alpha_2\nonumber\hspace*{-35pt}
\end{eqnarray}
for any $\alpha_1 \geq0$ and $\alpha_2 \geq0$ such that $0 \leq
\alpha_1 + \alpha_2 < 1$.
\end{theorem}

Theorem~\ref{theorembootmean} generalizes in the same way
that Theorem~\ref{theoremsubmean} generalizes. In particular, we have
the following result:
%
\begin{theorem} \label{theorembootgenmean}
Let $\mathbf P$ be defined as in Theorem~\ref{theoremsubmean}. Let
$J_n(x,P)$ be the distribution of the root (\ref{eqgeneralf}).
Suppose $f$ is continuous. Suppose further that for all $x \in\mathbf R$
%
\begin{eqnarray}
\label{eqonemoreone} P_n\bigl\{f\bigl(Z_n(P_n),
\Omega(\hat P_n)\bigr) \leq x\bigr\} &\rightarrow& P\bigl\{f(Z,
\Omega) \leq x\bigr\},
\\
\label{eqonemoretwo} P_n\bigl\{f\bigl(Z_n(P_n),
\Omega(\hat P_n)\bigr) < x\bigr\} &\rightarrow& P\bigl\{f(Z,\Omega) <
x\bigr\}
\end{eqnarray}
for any sequence $\{P_n \in\mathbf P\dvtx  n \geq1\}$ such that $Z_n(P_n)
\stackrel{d}{\rightarrow} Z$ under $P_n$ and $\Omega(\hat P_n)
\stackrel{P_n}{\rightarrow} \Omega$, where $Z \sim N(0,\Omega)$. Then
%
\begin{eqnarray}
\label{eqgenbootmeancover}
&&
\liminf_{n \rightarrow\infty} \inf_{P \in
\mathbf P} P\bigl\{
J_n^{-1}(\alpha_1,\hat P_n) \leq f
\bigl(Z_n(P),\hat\Omega_n\bigr) \leq
J_n^{-1}(1 - \alpha_2,\hat P_n)
\bigr\}\nonumber\\[-8pt]\\[-8pt]
&&\qquad
\geq1 - \alpha_1 - \alpha_2\nonumber
\end{eqnarray}
for any $\alpha_1 \geq0$ and $\alpha_2 \geq0$ such that $0 \leq
\alpha_1 + \alpha_2 < 1$.
\end{theorem}
\end{example}
%
\begin{example}[(Moment inequalities)]
Let $X^{(n)} = (X_1,\ldots, X_n)$ be an i.i.d. sequence of random
variables with distribution $P \in\mathbf P$ on $\mathbf R^k$ and
define $\mathbf P_0$ and $\mathbf P_1$ as in Example
\ref{exmomineqsub}. \citet{andrewsbarwick2012} propose testing the
null hypothesis that $P \in\mathbf P_0$ versus the alternative
hypothesis that $P \in\mathbf P_1$ at level $\alpha\in(0,1)$ using an
``adjusted quasi-likelihood ratio'' statistic $T_n(X^{(n)})$ defined as
follows:
\[
T_n\bigl(X^{(n)}\bigr) = \inf_{t \in\mathbf R^k\dvtx  t \leq0}
W_n(t)'\tilde\Omega_n^{-1}
W_n(t).
\]
Here, $t \leq0$ is understood to mean that the inequality holds component-wise,
\[
W_n(t) = \biggl( \frac{\sqrt n (\bar X_{1,n} - t_1)}{S_{1,n}},\ldots,
\frac{\sqrt n(\bar X_{k,n} - t_k)}{S_{k,n}}
\biggr)'
\]
and
%
\begin{equation}
\label{equationtildeomega} \tilde\Omega_n = \max\bigl\{\varepsilon-
\operatorname{det}(\hat\Omega_n), 0\bigr\} I_k + \hat
\Omega_n,
\end{equation}
where $\varepsilon> 0$ and $I_k$ is the $k$-dimensional identity matrix.
\citet{andrewsbarwick2012} propose a procedure for constructing
critical values for $T_n(X^{(n)})$ that they term ``refined moment
selection.'' For illustrative purposes, we instead consider in the
following theorem a simpler construction.
%
\begin{theorem} \label{theorembootmomineq}
Let $\mathbf P$ be defined as in Theorem~\ref{theoremsubmean}. Let
$J_n(x,P)$ be the distribution of the root
%
\begin{equation}
\label{eqaqlrroot} R_n\bigl(X^{(n)},P\bigr) =
\inf_{t \in\mathbf R^k\dvtx  t \leq0}\bigl(Z_n(P) - t\bigr)'\tilde
\Omega_n^{-1}\bigl(Z_n(P) - t\bigr),
\end{equation}
where $Z_n(P)$ is defined as in (\ref{eqznp}). Then, the test defined
by
\[
\phi_n\bigl(X^{(n)}\bigr) = I\bigl\{T_n
\bigl(X^{(n)}\bigr) > J_n^{-1}(1 - \alpha, \hat
P_n)\bigr\}
\]
satisfies (\ref{equniftest}) for any $0 < \alpha< 1$.
\end{theorem}

Theorem~\ref{theorembootmomineq} generalizes in a
straightforward fashion to other choices of test statistics, including
the one used in Theorem~\ref{theoremsubmomineq}. On the other hand,
even when the underlying choice of test statistic is the same, the
first-order asymptotic properties of the tests in Theorems \ref
{theorembootmomineq} and~\ref{theoremsubmomineq} will
differ. For other ways of constructing critical values that are more
similar to the construction given in \citet{andrewsbarwick2012}, see
\citet{romanoshaikhwolf2012}.
\end{example}
%
\begin{example}[(Multiple testing)]
Theorem~\ref{theoremboot} may be used in the same way that Theorem
\ref{theoremsubsample} was used in Example~\ref{exmultsub} to
construct tests of multiple hypotheses that behave well uniformly over
a large class of distributions.
To see this, let $X^{(n)} = (X_1,\ldots, X_n)$ be an i.i.d. sequence
of random variables with distribution $P \in\mathbf P$ on $\mathbf
R^k$, and again consider testing the family of null hypotheses (\ref
{eqfamily}) versus the alternative hypotheses (\ref{eqaltfamily}) in
a way that satisfies (\ref{equniffwer}) for $\alpha\in(0,1)$. For
$K \subseteq\{1,\ldots, k\}$, let $J_n(x,K,P)$ be the distribution
of the root
\[
R_n\bigl(X^{(n)},P\bigr) = \max_{j \in K}
\frac{\sqrt n (\bar X_{j,n} - \mu_j(P))}{S_{j,n}}
\]
under $P$, and consider the stepwise multiple testing procedure given
by Algorithm~\ref{algstepsub} with $L_n^{-1}(1 - \alpha, K_j)$
replaced by $J_n^{-1}(1 - \alpha, K_j, \hat P_n)$. We have the
following theorem, which is a bootstrap counterpart to Theorem \ref
{theoremsubmultiple}:
%
\begin{theorem} \label{theorembootmultiple}
Let $\mathbf P$ be defined as in Theorem~\ref{theoremsubmean}. Then
Algorithm~\ref{algstepsub} with $L_n^{-1}(1 - \alpha, K_j)$ replaced
by $J_n^{-1}(1 - \alpha, K_j, \hat P_n)$ satisfies (\ref
{eqfwercontrol}) for any $0 < \alpha< 1$.
\end{theorem}

It is, of course, possible to extend the analysis in a
straightforward way to two-sided testing.
\end{example}
%
\begin{example}[(Empirical process on $\mathbf R$)]
Let $X^{(n)} = (X_1,\ldots, X_n)$ be an i.i.d. sequence of random
variables with distribution $P \in\mathbf P$ on $\mathbf R$. Suppose
one wishes to construct a confidence region for the cumulative
distribution function associated with $P$, that is, $P\{(-\infty,t]\}
$. As described in Example~\ref{exsubemp}, a natural choice of root
in this case is given by (\ref{eqks}). In this setting, we have the
following theorem, which is a bootstrap counterpart to Theorem~\ref
{theoremsubemp}:
%
\begin{theorem} \label{theorembootemp}
Fix any $\varepsilon\in(0,1)$, and let $\mathbf P$ be defined as in
Theorem~\ref{theoremsubemp}. Let $J_n(x,P)$ be the distribution of
the root (\ref{eqks}). Denote by $\hat P_n$ the empirical
distribution of $X^{(n)}$. Then
\begin{eqnarray*}
&&\lim_{n \rightarrow\infty} \inf_{P \in\mathbf{P}} P \Bigl\{ J_n^{-1}(
\alpha_1, \hat P_n) \leq\sup_{t \in\mathbf R} \sqrt n \bigl|\hat
P_n\bigl\{(-\infty,t]\bigr\} - P\bigl\{(-\infty,t]\bigr\}\bigr|\\
&&\qquad\hspace*{182pt} \leq
J_n^{-1}(1-\alpha_2, \hat P_n)
\Bigr\}\\
&&\qquad = 1 - \alpha_1 - \alpha_2
\end{eqnarray*}
for any $\alpha_1 \geq0$ and $\alpha_2 \geq0$ such that $0 \leq
\alpha_1 + \alpha_2 < 1$.
\end{theorem}

Some of the conclusions of Theorem~\ref{theorembootemp} can
be found in \citet{romano1989unif}, though the method of proof given
in \citet{romanoshaikh2012supp} is quite different.
\end{example}
%
\begin{example}[(One sample $U$-statistics)] \label{exbootustat}
Let $X^{(n)} = (X_1,\ldots, X_n)$ be an i.i.d. sequence of random
variables with distribution $P \in\mathbf P$ on $\mathbf R$ and let
$h$ be a symmetric kernel of degree $m$. Suppose one wishes to
construct a confidence region for $\theta(P) = \theta_h(P)$ given by
(\ref{eqthetaustat}). As described in Example~\ref{exsubustat}, a
natural choice of root in this case is given by (\ref{equstatroot}).
Before proceeding, it is useful to introduce the following notation.
For an arbitrary kernel $\tilde h$, $\varepsilon> 0$ and $B > 0$, denote
by $\mathbf P_{\tilde h,\varepsilon,B}$ the set of all distributions $P$
on $\mathbf R$ such that
%
\begin{equation}
\label{eqcond1gen} E_P\bigl[\bigl|\tilde h(X_1,\ldots,X_m) - \theta_{\tilde h}(P)\bigr|^\varepsilon\bigr] \leq B.
\end{equation}
Similarly, for an arbitrary kernel $\tilde h$ and $\delta> 0$, denote
by $\mathbf S_{\tilde h,\delta}$ the set of all distributions $P$ on
$\mathbf R$ such that
%
\begin{equation}
\label{eqcond2gen} \sigma^2_{\tilde h}(P) \geq\delta,
\end{equation}
where $\sigma^2_{\tilde h}(P)$ is defined as in (\ref{eqsigmaP}).
Finally, for an arbitrary kernel $\tilde h$, $\varepsilon> 0$ and $B >
0$, let $\bar{\mathbf P}_{\tilde h,\varepsilon,B}$ be the set of
distributions $P$ on $\mathbf R$ such that
\[
E_P\bigl[\bigl|\tilde h(X_{i_1},\ldots,X_{i_m})
- \theta_{\tilde
h}(P)\bigr|^{\varepsilon}\bigr] \leq B,
\]
whenever $1 \leq i_j \leq n$ for all $1 \leq j \leq m$. Using this
notation, we have the following theorem:
%
\begin{theorem} \label{theorembootustat}
Define the kernel $h'$ of degree $2m$ according to the rule
%
\begin{eqnarray}
\label{eqhprime} h'(x_1,\ldots,x_{2m}) &=&
h(x_1,\ldots, x_m)h(x_1,x_{m+2},\ldots,x_{2m})\nonumber\\[-8pt]\\[-8pt]
&&{} - h(x_1,\ldots, x_m)h(x_{m+1},\ldots,x_{2m}).\nonumber
\end{eqnarray}
Suppose
\[
\mathbf P \subseteq\mathbf P_{h,2+\delta,B} \cap\mathbf S_{h,\delta
} \cap
\bar{ \mathbf P}_{h',1+\delta,B} \cap\bar{\mathbf P}_{h,2 +
\delta,B}
\]
for some $\delta> 0$ and $B > 0$. Let $J_n(x,P)$ be the distribution
of the root $R_n$ defined by (\ref{equstatroot}). Then
\[
\lim_{n \rightarrow\infty} \inf_{P \in\mathbf{P}} P \bigl\{ J_n^{-1}(
\alpha_1, \hat P_n) \leq\sqrt n\bigl(\hat
\theta_n - \theta(P)\bigr) \leq J_n^{-1}(1-
\alpha_2, \hat P_n) \bigr\} = 1 - \alpha_1 -
\alpha_2
\]
for any $\alpha_1$ and $\alpha_2$ such that $0 \leq\alpha_1 +
\alpha_2 < 1$.
\end{theorem}

Note that the kernel $h'$ defined in (\ref{eqhprime})
arises in the analysis of the estimated variance of the $U$-statistic.
Note further that the conditions on $\mathbf P$ in Theorem \ref
{theorembootustat} are stronger than the conditions on $\mathbf P$ in
Theorem~\ref{theoremsubustat}. While it may be possible to weaken the
restrictions on $\mathbf P$ in Theorem~\ref{theorembootustat} some,
it is not possible to establish the conclusions of Theorem \ref
{theorembootustat} under the conditions on $\mathbf P$ in Theorem~\ref
{theoremsubustat}. Indeed, as shown by \citet{bickelfreedman1981},
the bootstrap based on the root $R_n$ defined by (\ref{equstatroot})
need not be even pointwise asymptotically valid under the conditions on
$\mathbf P$ in Theorem~\ref{theoremsubustat}.
\end{example}

\begin{appendix}\label{app}
\section*{Appendix}

\subsection{\texorpdfstring{Proof of Theorem \protect\ref{theoremsubsample}}{Proof of Theorem 2.1}}
%
\begin{lemma} \label{lemmaquant}
If $F$ and $G$ are (nonrandom) distribution functions on $\mathbf{R}$,
then we have that:
\begin{longlist}[(ii)]
\item[(i)] If $\sup_{x \in\mathbf{R}} \{G(x) - F(x)\} \leq\varepsilon
$, then $G^{-1}(1 - \alpha_2) \geq F^{-1}(1 - (\alpha_2 + \varepsilon))$.
\item[(ii)] If $\sup_{x \in\mathbf{R}} \{F(x) - G(x)\} \leq
\varepsilon$, then $G^{-1}(\alpha_1) \leq F^{-1}(\alpha_1 + \varepsilon)$.
\end{longlist}
Furthermore, if $X \sim F$, it follows that:
\begin{longlist}[(iii)]
\item[(iii)] If $\sup_{x \in\mathbf{R}} \{G(x) - F(x)\} \leq
\varepsilon$, then $P\{X \leq G^{-1}(1 - \alpha_2)\} \geq1 - (\alpha_2
+ \varepsilon)$.
\item[(iv)] If $\sup_{x \in\mathbf{R}} \{F(x) - G(x)\} \leq
\varepsilon$, then $P\{X \geq G^{-1}(\alpha_1)\} \geq1 - (\alpha_1 +
\varepsilon)$.
\item[(v)] If $\sup_{x \in\mathbf{R}} |G(x) - F(x)| \leq\frac
{\varepsilon}{2}$, then $P\{G^{-1}(\alpha_1) \leq X \leq G^{-1}( 1 -
\alpha_2)\} \geq1 - (\alpha_1 + \alpha_2 + \varepsilon)$.
\end{longlist}
If $\hat G$ is a random distribution function on $\mathbf
{R}$, then we have further that:
\begin{longlist}[(viii)]
\item[(vi)] If $P\{\sup_{x \in\mathbf{R}} \{\hat G(x) - F(x)\} \leq
\varepsilon\} \geq1 - \delta$, then $P\{ X \leq\hat G^{-1}( 1 - \alpha
_2)\} \geq1 - (\alpha_2 + \varepsilon+
\delta)$.
\item[(vii)] If $P\{\sup_{x \in\mathbf{R}} \{F(x) - \hat G(x)\}
\leq\varepsilon\} \geq1 - \delta$, then $P\{ X \geq\hat G^{-1}(\alpha_1)
\} \geq1 - (\alpha_1 + \varepsilon+
\delta)$.
\item[(viii)] If $P\{\sup_{x \in\mathbf{R}} |\hat G(x) - F(x)| \leq
\frac{\varepsilon}{2}\} \geq1 - \delta$, then $P\{\hat G^{-1}(\alpha_1)
\leq X \leq\hat G^{-1}( 1 - \alpha_2)\} \geq1 - (\alpha_1 +
\alpha_2 + \varepsilon+ \delta)$.
\end{longlist}
\end{lemma}
\begin{pf}
To see (i), first note that $\sup_{x \in
\mathbf{R}} \{G(x) - F(x)\} \leq\varepsilon$ implies that $G(x) -
\varepsilon\leq F(x)$ for all $x \in\mathbf{R}$. Thus, $\{x \in\mathbf
R\dvtx  G(x) \geq1 - \alpha_2 \} = \{x \in\mathbf R\dvtx  G(x) - \varepsilon
\geq1 - \alpha_2 - \varepsilon\} \subseteq\{x \in\mathbf R\dvtx  F(x)
\geq1 - \alpha_2 - \varepsilon\}$, from which it follows that
$F^{-1}(1- (\alpha_2 + \varepsilon)) = \inf\{x \in\mathbf R\dvtx  F(x)
\geq1 - \alpha_2 - \varepsilon\} \leq\inf\{x \in\mathbf R\dvtx  G(x)
\geq1 - \alpha_2\} = G^{-1}(1 - \alpha_2 )$. Similarly, to prove
(ii), first note that $\sup_{x \in\mathbf{R}} \{F(x) - G(x)\} \leq
\varepsilon$ implies that $F(x) - \varepsilon\leq G(x)$ for all $x \in
\mathbf{R}$, so $\{x \in\mathbf R\dvtx  F(x) \geq\alpha_1 + \varepsilon\}
= \{x \in\mathbf R\dvtx  F(x) - \varepsilon\geq\alpha_1\} \subseteq\{x\in
\mathbf R\dvtx  G(x) \geq\alpha_1\}$. Therefore, $G^{-1}(\alpha_1) =
\inf\{ x \in\mathbf R\dvtx  G(x) \geq\alpha_1\} \leq\inf\{x \in
\mathbf R\dvtx  F(x) \geq\alpha_1 + \varepsilon\} = F^{-1}(\alpha_1 +
\varepsilon)$. To prove (iii), note that because $\sup_{x \in\mathbf
{R}} \{G(x) - F(x)\} \leq\varepsilon$, it follows from (i) that $\{X
\leq G^{-1}( 1 - \alpha_2)\} \supseteq\{X \leq F^{-1}( 1 - (\alpha_2
+ \varepsilon))\}$. Hence, $P\{X \leq G^{-1}( 1 - \alpha_2)\} \geq P\{X
\leq F^{-1}(1 - (\alpha_2 + \varepsilon))\} \geq1 - (\alpha_2 +
\varepsilon)$. Using the same reasoning, (iv) follows from (ii) and the
assumption that $\sup_{x \in\mathbf{R}} \{F(x) - G(x)\} \leq
\varepsilon$. To see (v), note that
\begin{eqnarray*}
P\bigl\{G^{-1}(\alpha_1) \leq X \leq G^{-1}(1 -
\alpha_2)\bigr\} &\geq& 1 - P\bigl\{ X < G^{-1}(
\alpha_1)\bigr\} \\
&&{}- P\bigl\{X > G^{-1}(1 -
\alpha_2)\bigr\}
\\
&\geq& 1 - (\alpha_1 + \alpha_2
+ \varepsilon),
\end{eqnarray*}
where the first inequality follows from the Bonferroni inequality, and
the second inequality follows from (iii) and (iv). To prove (vi), note that
\begin{eqnarray*}
&&
P\bigl\{X \leq\hat G^{-1}(1 - \alpha_2)\bigr\} \\
&&\qquad\geq P
\Bigl\{X \leq\hat G^{-1}(1 - \alpha_2) \cap
\sup_{x \in\mathbf{R}} \bigl\{ \hat G(x) - F(x)\bigr\} \leq\varepsilon
\Bigr\}
\\
&&\qquad\geq P \Bigl\{X \leq F^{-1}\bigl(1 - (\alpha_2 +
\varepsilon)\bigr) \cap\sup_{x \in\mathbf{R}} \bigl\{ \hat G(x) - F(x)\bigr
\} \leq
\varepsilon\Bigr\}
\\
&&\qquad\geq P \bigl\{X \leq F^{-1}
\bigl(1 - (\alpha_2 + \varepsilon)\bigr)\bigr\} - P\Bigl\{
\sup_{x \in\mathbf{R}} \bigl\{ \hat G(x) - F(x)\bigr\} > \varepsilon\Bigr
\}
\\
&&\qquad= 1 - \alpha_2 -\varepsilon-\delta,
\end{eqnarray*}
where the second inequality follows from (i). A similar argument using
(ii) establishes (vii). Finally, (viii) follows from (vi) and (vii) by
an argument analogous to the one used to establish (v).
\end{pf}
%
\begin{lemma} \label{lemmavcreal}
Let $X^{(n)} = (X_1,\ldots, X_n)$ be an i.i.d. sequence of random
variables with distribution $P$. Denote by $J_n(x,P)$ the distribution
of a real-valued root $R_n = R_n(X^{(n)},P)$ under $P$. Let $N_n = {n
\choose b}$, $k_n = \lfloor\frac{n}{b} \rfloor$ and define
$L_n(x,P)$ according to (\ref{equationsubdist}). Then, for any
$\varepsilon> 0$, we have that
%
\begin{equation}
\label{equationbound2} P \Bigl\{\sup_{x \in\mathbf{R}} \bigl|L_n(x,P) -
J_b(x,P)\bigr| > \varepsilon\Bigr\} \leq\frac{1}{\varepsilon}\sqrt{
\frac{2 \pi}{k_n}}.
\end{equation}
%
\end{lemma}
\begin{pf}
Let $\varepsilon> 0$ be given and define $S_n(x,
P; X_1,\ldots, X_n)$ by
\[
\frac{1}{k_n} \sum_{1 \leq i \leq k_n} I\bigl
\{R_b\bigl((X_{b(i - 1) + 1},\ldots, X_{bi}),P\bigr) \leq
x\bigr\} - J_b(x,P).
\]
Denote by $\mathcal{S}_n$ the symmetric group with $n$ elements. Note
that using this notation, we may rewrite $L_n(x,P) - J_b(x,P)$ as
\[
Z_n(x,P; X_1,\ldots, X_n) =
\frac{1}{n!} \sum_{\pi\in\mathcal
{S}_n} S_n(x,P;
X_{\pi(1)},\ldots, X_{\pi(n)}).
\]
Note further that
\[
\sup_{x \in\mathbf{R}} \bigl|Z_n(x,P; X_1,\ldots,
X_n)\bigr| \leq\frac
{1}{n!} \sum_{\pi\in\mathcal{S}_n}
\sup_{x \in\mathbf{R}} \bigl|S_n(x,P; X_{\pi(1)},\ldots,
X_{\pi(n)})\bigr|,
\]
which is a sum of $n!$ identically distributed random variables. Let
$\varepsilon> 0$ be given. It follows that $P\{ \sup_{x \in\mathbf{R}}
|Z_n(x,P; X_1,\ldots, X_n)| > \varepsilon\}$ is bounded above by
%
\begin{equation}
\label{equationdumb} P \biggl\{ \frac{1}{n!} \sum
_{\pi\in\mathcal{S}_n} \sup_{x \in
\mathbf{R}} \bigl|S_n(x,P;
X_{\pi(1)},\ldots, X_{\pi(n)})\bigr| > \varepsilon\biggr\}.
\end{equation}
Using Markov's inequality, (\ref{equationdumb}) can be bounded by
%
\begin{eqnarray}
\label{equationintreal}
&&
\frac{1}{\varepsilon}E_P \Bigl[
\sup_{x \in\mathbf{R}} \bigl|S_n(x,P; X_1,\ldots,
X_n)\bigr| \Bigr]\nonumber\\[-8pt]\\[-8pt]
&&\qquad = \frac{1}{\varepsilon} \int_0^1
P \Bigl\{ \sup_{x \in\mathbf{R}} \bigl|S_n(x,P; X_1,\ldots,
X_n)\bigr| > u \Bigr\} \,du.\nonumber
\end{eqnarray}
We may use the Dvoretsky--Kiefer--Wolfowitz inequality to bound the
right-hand side of (\ref{equationintreal}) by
\[
\frac{1}{\varepsilon}\int_0^1 2 \exp\bigl\{- 2
k_n u^2\bigr\} \,du = \frac
{2}{\varepsilon}\sqrt{
\frac{2 \pi}{k_n}} \biggl[\Phi(2\sqrt{k_n}) - \frac{1}{2}
\biggr] < \frac{1}{\varepsilon}\sqrt{\frac{2 \pi}{k_n}},
\]
which establishes (\ref{equationbound2}).
\end{pf}


%
\begin{lemma} \label{lemmafinite}
Let $X^{(n)} = (X_1,\ldots, X_n)$ be an i.i.d. sequence of random
variables with distribution $P \in\mathbf P$. Denote by $J_n(x,P)$ the
distribution of a real-valued\vadjust{\goodbreak} root $R_n = R_n(X^{(n)},P)$ under $P$.
Let $k_n = \lfloor\frac{n}{b} \rfloor$ and define $L_n(x,P)$
according to~(\ref{equationsubdist}). Let
\begin{eqnarray*}
\delta_{1,n}(\varepsilon,\gamma,P) &=& \frac{1}{\gamma\varepsilon}\sqrt{
\frac{2 \pi}{k_n}} + I \Bigl\{ \sup_{x \in\mathbf{R}} \bigl\{J_b(x,P) -
J_n(x,P)\bigr\} > (1- \gamma)\varepsilon\Bigr\},
\\
\delta_{2,n}(\varepsilon,\gamma,P) &=& \frac{1}{\gamma\varepsilon}\sqrt{
\frac{2 \pi}{k_n}} + I \Bigl\{ \sup_{x \in\mathbf{R}} \bigl\{J_n(x,P) -
J_b(x,P)\bigr\} > (1- \gamma)\varepsilon\Bigr\},
\\
\delta_{3,n}(\varepsilon,\gamma,P) &=& \frac{1}{\gamma\varepsilon}\sqrt{
\frac{2 \pi}{k_n}} + I \Bigl\{\sup_{x \in\mathbf{R}} \bigl|J_b(x,P) -
J_n(x,P)\bigr| > (1- \gamma)\varepsilon\Bigr\}.
\end{eqnarray*}
Then, for any $\varepsilon> 0$ and $\gamma\in(0,1)$, we have that:
\begin{longlist}[(iii)]
\item[(i)] $P\{R_n \leq L_n^{-1}(1 - \alpha_2,P)\} \geq1 - (\alpha_2 +
\varepsilon+ \delta_{1,n}(\varepsilon, \gamma,
P))$;
\item[(ii)] $P\{R_n \geq L_n^{-1}(\alpha,P)\} \geq1 - (\alpha_1 +
\varepsilon+ \delta_{2,n}(\varepsilon, \gamma, P))$;\vspace*{1pt}
\item[(iii)] $P\{L_n^{-1}(\alpha_1,P) \leq R_n \leq L_n^{-1}(1 -
\alpha_2,P)\} \geq1 - (\alpha_1 + \alpha_2 + \varepsilon+ \delta
_{3,n}(\varepsilon, \gamma,
P))$.
\end{longlist}
\end{lemma}
\begin{pf}
Let $\varepsilon> 0$ and $\gamma\in(0,1)$ be
given. Note that
\begin{eqnarray*}
& & P \Bigl\{ \sup_{x \in\mathbf{R}} \bigl\{L_n(x,P) -
J_n(x,P)\bigr\} > \varepsilon\Bigr\}
\\
&&\qquad\leq P \Bigl\{\sup_{x \in\mathbf{R}} \bigl\{L_n(x,P) -
J_b(x,P)\bigr\} + \sup_{x \in\mathbf{R}} \bigl\{J_b(x,P) -
J_n(x,P)\bigr\} > \varepsilon\Bigr\}
\\
&&\qquad\leq P \Bigl\{\sup_{x \in\mathbf{R}} \bigl\{L_n(x,P) -
J_b(x,P)\bigr\} > \gamma\varepsilon\Bigr\} \\
&&\qquad\quad{}+ I \Bigl\{
\sup_{x \in\mathbf{R}} \bigl\{ J_b(x,P) - J_n(x,P)\bigr\} >
(1- \gamma)\varepsilon\Bigr\}
\\
&&\qquad\leq \frac{1}{\gamma\varepsilon}\sqrt{\frac{2 \pi}{k_n}} + I \Bigl\{
\sup_{x \in\mathbf{R}} \bigl\{J_b(x,P) - J_n(x,P)\bigr\} >
(1- \gamma)\varepsilon\Bigr\},
\end{eqnarray*}
where the final inequality follows from Lemma~\ref{lemmavcreal}.
Assertion (i) thus follows from the definition of $\delta_{1,n}(\varepsilon
,\gamma,P)$ and part (vi) of Lemma~\ref{lemmaquant}.
Assertions~(ii) and (iii) are established similarly.
\end{pf}
\begin{pf*}{Proof of Theorem~\ref{theoremsubsample}}
To prove (i), note that by part (i) of Lem\-ma~\ref{lemmafinite}, we have
for any $\varepsilon> 0$ and $\gamma\in(0,1)$ that
\[
\sup_{P \in\mathbf P}P\bigl\{R_n \leq L_n^{-1}(1
- \alpha_2, P)\bigr\} \geq1 - \Bigl(\alpha_2 + \varepsilon+
\inf_{P \in\mathbf P} \delta_{1,n}(\varepsilon, \gamma, P) \Bigr),
\]
where
\[
\delta_{1,n}(\varepsilon,\gamma,P) = \frac{1}{\gamma\varepsilon}\sqrt{
\frac{2 \pi}{k_n}} + I \Bigl\{ \sup_{x \in\mathbf{R}} \bigl\{J_b(x,P) -
J_n(x,P)\bigr\} > (1- \gamma)\varepsilon\Bigr\}.
\]
By the assumption on $\sup_{P \in\mathbf{P}} \sup_{x \in\mathbf
{R}} \{J_b(x,P) - J_n(x,P) \}$, we have that $\inf_{P \in\mathbf P}
\delta_{1,n}(\varepsilon,\gamma,P) \rightarrow0$ for every $\varepsilon>
0$. Thus, there exists a sequence $\varepsilon_n > 0$ tending to $0$ so
that $\inf_{P \in\mathbf P} \delta_{1,n}(\varepsilon_n,\gamma,P)
\rightarrow0$. The desired claim now follows from applying part (i) of
Lemma~\ref{lemmafinite} to this sequence. Assertions (ii) and~(iii)
follow in exactly the same way.
\end{pf*}

\subsection{\texorpdfstring{Proof of Theorem \protect\ref{theoremboot}}{Proof of Theorem 2.4}}

We prove only (i). Similar arguments can be used to establish (ii) and
(iii). Let $\alpha_1 = 0$, $0 \leq\alpha_2 < 1$ and $\eta> 0$ be
given. Choose $\delta> 0$ so that
\[
\sup_{x \in\mathbf{R}} \bigl\{J_n\bigl(x,P'\bigr) -
J_n(x,P)\bigr\} < \frac{\eta}{2},
\]
whenever $\rho(P',P) < \delta$ for $P' \in\mathbf P'$ and $P \in
\mathbf P$. For $n$ sufficiently large, we have that
\[
\sup_{P \in\mathbf{P}} P\bigl\{\rho(\hat P_n, P) > \delta\bigr\} <
\frac
{\eta}{4} \quad\mbox{and}\quad \sup_{P \in\mathbf{P}} P\bigl\{\hat P_n
\notin\mathbf P'\bigr\} < \frac{\eta}{4}.
\]
For such $n$, we therefore have that
\begin{eqnarray*}
1 - \frac{\eta}{2} &\leq&\inf_{P \in\mathbf{P}} P\bigl\{\rho(\hat
P_n,P) \leq\delta\cap\hat P_n \in\mathbf
P'\bigr\} \\
&\leq&\inf_{P \in
\mathbf{P}} P \biggl\{\sup_{x \in\mathbf{R}}
\bigl\{J_n(x,\hat P_n) - J_n(x,P)\bigr\}
\leq\frac{\eta}{2} \biggr\}.
\end{eqnarray*}
It follows from part (vi) of Lemma~\ref{lemmaquant} that for such $n$
\[
\inf_{P \in\mathbf{P}} P\bigl\{R_n \leq J_n^{-1}(1
- \alpha_2,\hat P_n)\bigr\} \geq1 - (
\alpha_2 + \eta).
\]
Since the choice of $\eta$ was arbitrary, the desired result follows.
\end{appendix}

\begin{supplement}
\stitle{Supplement to ``On the uniform asymptotic validity of subsampling and the bootstrap''}
\slink[doi]{10.1214/12-AOS1051SUPP} 
\sdatatype{.pdf}
\sfilename{aos1051\_supp.pdf}
\sdescription{The supplement provides additional details and proofs for
many of the results in the authors' paper.}
\end{supplement}


\printaddresses

\end{document}